%% file: main_onecol.latex
\documentclass[10pt]{article}
\usepackage[a4paper, total={6.5in,9in}]{geometry}
\usepackage{cite}     
\usepackage{graphicx} 
\usepackage{psfrag}
\usepackage{subfigure} 
\usepackage{url}       
\usepackage{amsmath}   
\usepackage{amsthm}   
\usepackage{algpseudocode}
\usepackage{enumerate}
\usepackage{amssymb}
\usepackage{units}
\usepackage[english]{babel}
\addto\captionsenglish{}

\DeclareMathOperator{\real}{Re}
\DeclareMathOperator{\imag}{Im}

\begin{document}
\title{\LARGE \bf Scalable Power System Line Upgrade Planning With Policy Constraints: A ``Branch and Benders'' Approach}
\author{ 
  Sandro Merkli\thanks{
    embotech AG, Giessereistrasse 18, 8005 Zurich}
  \thanks{
    Automatic Control Lab,
    ETH Zurich, Physikstrasse 3, 8092 Zurich
    \tt{ \{smerkli,rsmith\}@control.ee.ethz.ch}}\;,
  Roy S.\ Smith\footnotemark[2]
}

\include{definitions}
\date{Latest revision: January 2022}
\newtheorem{lemma}{Lemma}

\maketitle
\thispagestyle{empty}
\pagestyle{empty}
\input{content}
\bibliographystyle{ieeetr}
\bibliography{biblio}
\end{document}

%% file: definitions.tex
\def \todo {\textbf{todo:} }
\def \div { \mathrm{div} }
\def \rot { \mathrm{\textbf{rot}} }
\def \minfty { -\infty }
\def \infint { \int_{-\infty}^\infty }
\def \infsum { \sum_{k = -\infty}^\infty }
\def \adj { \mathrm{adj} }
\def \lps { \quad \laplace \quad }
\def \B { \mathbb B }
\def \C { \mathbb C }
\def \N { \mathbb N }
\def \R { \mathbb R }
\def \Z { \mathbb Z }
\def \Q { \mathbb Q }
\def \Ac { \mathcal A }
\def \Bc { \mathcal B }
\def \Dc { \mathcal D }
\def \Ec { \mathcal E }
\def \Fc { \mathcal F }
\def \Gc { \mathcal G }
\def \Xc { \mathcal X }
\def \Ic { \mathcal I }
\def \Zc { \mathcal Z }
\def \Qc { \mathcal Q }
\def \Rc { \mathcal R }
\def \Uc { \mathcal U }
\def \Nc { \mathcal N }
\def \ei { \varepsilon_{\text{int}} }
\def \xtt { x_{t|t}   }
\def \gdw { \;\; \Longleftrightarrow \;\; }
\def \diag { \operatorname{diag} }
\def \minim {\operatorname*{minimize} }
\def \maxim {\operatorname*{maximize} }
\def \st { \operatorname*{subject\ to} }
\def \dom { \operatorname{dom} }
\def \img { \operatorname{im} }
\def \conv { \operatorname{conv} }
\def \cone { \operatorname{cone} }
\def \projop { \operatorname{proj} }
\def \jw  { j\omega }
\def \bmb { \begin{bmatrix} }
\def \bme { \end{bmatrix} }
\newcommand{\cve}[1] {#1}
\newcommand{\jcve}[1] {\bar{#1}}
\newcommand{\drop}[1] {{ }}

%% file: content.tex
\begin{abstract}
  The integration of more renewable energy sources into the power system is
  presenting system operators with various challenges. At the distribution system
  level, voltage magnitudes that violate operating limits near large
  photovoltaic installations have been observed. While these issues can be
  partially mitigated with more advanced control, hardware upgrades are
  required at some point. This work presents a scalable, optimization-based
  approach for deciding which lines in a network to upgrade. Compared to
  existing approaches, it explicitly takes the operating policy of the system
  into account and provides both reasonable solutions in short computation times
  as well as globally optimal solutions when run to completion. Compared to
  earlier work on the same topic, an extended computational approach is taken
  that can simultaneously optimize for many load scenarios across arbitrary 
  configurations of machines and CPU cores per machine in a scalable manner by 
  using the Benders decomposition. In addition to the theory, numerical 
  experiments are presented along with a discussion of the scaling properties 
  of the Benders-based approach, giving potential users a better basis to decide 
  whether their problem is big enough for the approach to make sense. 
\end{abstract}

\section{Introduction}
The share of renewable energy sources is increasing in power
systems~\cite{ren212015global}. While newly deployed renewables make power
generation more sustainable, they also present challenges to power system
operators. One such challenge is voltage magnitude rise due to
photovoltaic systems in low to medium voltage distribution
networks~\cite{Ayres2010}. Because the lines in such networks have
non-negligible resistance, voltage differences between buses become relevant.
If the power in-feed due to local generation is too high, voltage magnitudes
can exceed the operating limits near the feed-in point.  In order to avoid damage
to end user devices, steps have to be taken to mitigate such overvoltage
occurrences.

System operators can either change the operating policy of the system or
perform upgrades to the system itself. Examples of the former include
curtailment of the photovoltaic systems~\cite{Merkli2018} and changes in
on-load tap changing (OLTC) transformer control schemes. Such changes tend to
be less costly to implement than system hardware upgrades. On the other hand,
curtailment is not desirable as valuable renewable energy remains unused. System
hardware upgrade possibilities include energy storage devices, OLTC hardware
upgrades or power line upgrades. In this work, only power line upgrades are
considered since they are the most readily available at the time of this
writing. 

The general optimization problem of deciding which components of a power system
to upgrade is hard to solve: 
The question whether an admissible operating point exists for a
given AC power system and load pattern is already NP-hard to answer in
general~\cite{Lehmann2016}. For line upgrades in particular, an additional
difficulty stems from the fact that line hardware is typically available for
purchase with fixed increments in admittance. Even if line hardware was more 
available with user-specifiable admittance, the cost of upgrading a line is in
practice dominated by the construction work required for the upgrade and not by
the line admittance cost. In order to model this cost accurately, integrality
constraints have to be added to the problem, further increasing its difficulty.

These complicating factors make the power system line upgrade problem
intractable to solve in general. Existing methods therefore take steps to
simplify the problem in various ways. In terms of objective, most existing work
minimizes upgrade or upgrade plus operation cost for a given
scenario~\cite{Rider2007,Ramirez-Rosado1998}. One line of research then
approximates the more accurate AC model by a computationally tractable DC model
consisting of linear equations~\cite{Macedo2016}. The resulting mixed-integer
linear problem, while still NP-hard, can usually be solved in a reasonable
amount of time with state-of-the-art optimization solvers. The caveat in this
approach is that the result is not guaranteed to be a viable solution
to the full AC model of the system, let alone practically deployable. 
A related approach uses convex relaxations of the power flow equations, such
as the ones presented in~\cite{Low2014a,Gan2014}, to create formulations that
share the same caveat, but provide the added benefit of a lower bound on the 
number of line upgrades required~\cite{Jabr2013}.
A different approach is to apply heuristic methods for
finding reasonable solutions to the
problem~\cite{Ramirez-Rosado1998,Koutsoukis}. It has been shown that these
methods find reasonably good solutions to practical problems, but they do not
provide any lower bounds on the objective value.

The method presented in this paper combines the lower-bounding and heuristic
methods by means of a Branch-and-Bound procedure. This procedure is then
extended with constraint generation and a distributed optimization based on the
Benders decomposition. The result is an algorithm that can both find
practically deployable solutions to the line upgrade problem in reasonable
computation times as well as certified globally optimal solutions if it is run
to termination. The use of the Benders decomposition algorithm for solving 
relaxations enables the method to trivially be able to make use of a number of 
parallel computing agents up to the number of scenarios considered in the 
optimization. Since the most computationally demanding parts of the 
algorithm can be run in parallel this way, computation times are almost linearly
improved by the addition of more parallel agents. 

This scalability in the number of scenarios is a big step towards practical
applicability of sampling-based
methods~\cite{Esfahani2015,Campi2018,Calafiore2006} for power system planning
problems.

\subsection{Contribution}
The method presented improves on the state of the art in the following aspects:
\begin{enumerate}[a)]
  \item The operating policy used in practice is integrated directly into the
    optimization problem formulation. This guarantees that the selected upgrade
    configuration is also feasible when deployed, while existing approaches
    at best guarantee that a feasible operating point exists.
  \item The algorithm is compatible with existing heuristics for finding useful
    solutions (upper bounds) and can make use of existing convex relaxations.
  \item Multiple scenarios can be specified instead of a single worst-case
    scenario and the system topology can be optimized over all of them jointly.
\end{enumerate}
Compared to our earlier conference paper on the same
topic~\cite{merkli2018globally}, this work includes the following extensions:
\begin{enumerate}[1)]
  \item A detailed treatment of the reformulation of the original problem into
    a mixed-integer quadratically constrained quadratic problem (MI-QCQP) is
    given. 
  \item An approach based on the Benders decomposition is applied to the
    mixed-integer semidefinite (MI-SDP) relaxations arising in the
    Branch-and-Bound procedure and integrated with the latter, making the
    procedure scalable and amenable to parallelization. 
  \item Numerical experiments are presented that demonstrate how the described
    approach benefits from parallelism.
\end{enumerate}

\subsection{Contents}
Section~\ref{sec:modeling} outlines the model used in this work.
Section~\ref{sec:probbigM} presents the problem formulation and its
reformulation into a mixed-integer quadratic problem with policy constraints.
Section~\ref{sec:algorithm} then presents the main algorithm. Numerical
experiments are presented and discussed in Section~\ref{sec:numexp}.

\section{Modeling}
\label{sec:modeling}
In the following, $\bar x$ will be used to denote the element-wise complex
conjugate of a variable $x$. The operators $\real(x)$ and $\imag(x)$ are 
used to denote the real and imaginary parts of variable $x$. The operator
$x\cdot y$ is used in some places to denote the product $xy$ of scalars~$x$
and~$y$ to improve readability.

\subsection{Power system}
A power system with $N$ buses and $L$ lines is modeled using the bus 
injection model of the AC power flow equations,
\begin{equation}
  \label{eqn:kirchhoff}
  s = \diag(v)\overline{Yv},
\end{equation}
where $s \in \C^N$ and $v \in \C^N$ represent the complex powers and voltages
at the buses of the system. The Laplacian~$Y \in \C^{N\times N}$ contains all
information about the system topology considered in this work: 
\begin{equation}
  \label{eqn:laplacian}
  Y_{jl} := \begin{cases} 
    y_{jl} & \text{if } j \ne l, \\
    - \sum_{k=1,k \ne j}^N y_{jk} & \text{if } j = l,
  \end{cases}
\end{equation}
where $y_{jl}$ is the complex admittance of the line between buses $j$ and $l$.
Shunt admittances are not included in the model for notational simplicity, but
could be added by modifying
equations~\eqref{eqn:laplacian},~\eqref{eqn:opconst2} and the reformulation
in Section~\ref{sec:probbigM} accordingly.

\subsection{Operating limits}
All buses are assumed to have voltage magnitude limits as well as power limits:
\begin{equation}
  \label{eqn:opconst1}
  \begin{aligned}
    v_{\mathrm{min},j} &\le |v_j| \le v_{\mathrm{max},j},\\
    p_{\mathrm{min},j} &\le \real(s_j) \le p_{\mathrm{max},j}, \\
    q_{\mathrm{min},j} &\le \imag(s_j) \le q_{\mathrm{max},j},
  \end{aligned}
\end{equation}
By making different subsets of the above constraints tight (i.e. the minimum
value equal to the maximum value), one can model PV, PQ as well as slack buses
in a unified manner. For example, a PQ load could have the latter two limits
tight and the voltage magnitude constrained to $[0.9,1.1]$ per unit. For 
the slack bus, the voltage magnitude would be fixed but the other two limits
would not be present. The thermal limits on power lines are represented using a
current-based formulation:
\begin{equation}
  \label{eqn:opconst2}
  |Y_{jl}|\cdot|v_j-v_l| \le I_{\mathrm{max,jl}}.
\end{equation}
Note that other convex constraints, such as apparent power limits on
generators, can be included in the formulation if desired and the method
presented here is still applicable. 

\subsection{Violating scenarios}
Similar to other work in the field, we take the approach of minimizing the
cost of upgrades. However, instead of looking at only a worst case, we
introduce the concept of \emph{violating scenarios}. We assume the availability
of $K$ separate system steady-state scenarios, indexed with the letter $k$, in
the form of powers and voltages $(s^k,v^k)$ along with power limit data
$p^k_\text{min}$, $p^k_\text{max}$, $q^k_\text{min}$ and $q^k_\text{max}$.
These scenarios can come from past measurements, simulations of hypothetical
scenarios or also from worst-case studies. The method presented here is
scalable in the number of such scenarios. It will find a set of upgrades that,
combined with the operating policy introduced in the next section, leads to no
constraint violations in all scenarios. If no such set exists, the method will
certify that this is the case.

\subsection{Operating policies}
As mentioned in the introduction, it has been shown to be hard to decide
whether a power dispatch exists that satisfies the Kirchhoff
equations~\eqref{eqn:kirchhoff} as well as the operating
constraints~\eqref{eqn:opconst1}, \eqref{eqn:opconst2} for a given load
pattern. Fortunately, what is relevant in practice to a system operator is less 
whether a solution exists in theory, but rather whether they can find it
with their system operating policy. We model this policy as a function of the
system topology $Y$ and the operating limits specific to the scenario:
\begin{equation}
  \label{eqn:policy}
  (\tilde s^k,\tilde v^k) = g(Y,\text{limits}^k),
\end{equation}
where ``limits$^k$'' includes the data $v_\text{max}$, $v_\text{min}$, $I_\text{max}$,
$p_\text{min}^k$, $p_\text{max}^k$, $q_\text{min}^k$ and $q_\text{max}^k$. This work
assumes that the voltage and line current limits are the same across all
scenarios whereas the active and reactive power limits may be different.
The policy function can choose values for the powers and voltages within
the given limits, leading to new vectors $(\tilde s^k,\tilde v^k)$. This is
done to model the fact that, given a different topology, a different dispatch might
have been chosen. Examples of such policy functions are an AC economic dispatch
or local controllers. The assumptions made about $g$ are that it is
tractable to evaluate and that the $(\tilde s^k,\tilde v^k)$ returned by~\eqref{eqn:policy} 
satisfy~\eqref{eqn:kirchhoff}.

\subsection{Line upgrades}
Line upgrades are modeled as changes in admittance and can be formulated in 
a vectorized manner, 
\begin{equation}
  \label{eqn:yupg}
  Y_{\text{upg}}(a) = Y + \sum_{i=1}^{n_u} (a_i \cdot \Delta Y_i),
\end{equation}
where $a_i \in \{0,1\}$ indicates whether the upgrade is performed $(1)$ or not
$(0)$, the constant matrix $\Delta Y_i$ determines the change from the original
system topology~$Y$ if upgrade~$i$ is performed and $n_u$ is the total number 
of such upgrade possibilities. Constraints on the upgrade combinations are 
expressed as a separate set of linear constraints,
\begin{equation}
  \label{eqn:upgcon}
  Aa \le b.
\end{equation}
This set of constraints is used to ensure that at most one type of line upgrade
is performed per line. If a line upgrade is performed, its current limit can
also change,
\begin{equation}
  \label{eqn:imaxupg}
  I_{\max,\text{upg},jl} = I_{\max,jl} + \sum_{i \in \mathcal U_{jl}} 
    (a_i \cdot \Delta I_{jl,i}),
\end{equation}
where $\mathcal U_{jl}$ is the set of indexes $i$ which refer to upgrade
choices that affect the line from bus $j$ to bus $l$.

\section{Upgrade problem and Big-M formulation}
\label{sec:probbigM}
In this section, the upgrade problem is formulated mathematically. It is then
brought into a standard MI-QCQP form, which in turn will admit an MI-SDP
relaxation in Section~\ref{sec:algorithm}.

\subsection{Problem formulation and standard form} 
The system upgrade problem is given as follows:
\begin{subequations}
  \makeatletter
  \def\@currentlabel{U}
  \makeatother
  \label{eqn:upgradeproblem}
  \renewcommand{\theequation}{U.\arabic{equation}}
  \begin{align}
  \text{Problem U:} \nonumber \\ 
    \minim_{a,\tilde s^k,\tilde v^k} \;\;& 1^Ta \\
  \st \;\;& a \in \{0,1\}^{n_u}, \\
        & \eqref{eqn:policy}, \eqref{eqn:yupg}, 
          \eqref{eqn:upgcon}, \eqref{eqn:imaxupg}, \\
        & \eqref{eqn:kirchhoff}, \eqref{eqn:opconst1} \text{ for } 
          \tilde s^k,\tilde v^k, Y_\text{upg}(a), \\
        & \eqref{eqn:opconst2} \text{ for } \tilde v^k, I_{\max,\text{upg}}, \\
        &\;\; \forall k \in \{1,\ldots, K\} \nonumber
  \end{align}
\end{subequations}
where $a,\tilde v^k, \tilde s^k$ are optimization variables and
the remaining symbols are given data. The cost function is chosen to be the
sum of upgrades performed, but any convex cost function in $a$ would be
admissible. Problem~\eqref{eqn:upgradeproblem}, equivalently reformulated as an
MI-QCQP with an added policy constraint, is
\begin{subequations}
  \makeatletter
  \def\@currentlabel{P}
  \makeatother
  \label{eqn:upgradeproblem2}
  \renewcommand{\theequation}{P.\arabic{equation}}
  \begin{align}
  \text{Problem P:} \nonumber \\ 
  \minim_{a,z^k,y^k} &\;\; 1^Ta \\
         \st &\;\; a \in \{0,1\}^{n_u}\label{eqn:up2_bin}, \\ 
             &\;\; Aa \le b, \label{eqn:up2_poly}\\
             &\;\; Cy^k \le d^k, \label{eqn:up2_snappoly}\\
             &\;\; \alpha_h \le (z^k)^TQ_hz^k + q_h^Ty^k + m_h^Ta \le \beta_h,
               \label{eqn:up2_quad} \\
             &\;\; (\tilde s^k,\tilde v^k) = g(a,\text{limits}^k) \label{eqn:up2_policy} \\ 
             &\;\; \forall h \in \{1,\ldots, H\}, \forall k \in \{1,\ldots, K\}.
               \nonumber
  \end{align}
\end{subequations}
The following sections outline the variable correspondences and how the
constraints are brought into the standard form above.
Problem~\eqref{eqn:upgradeproblem2} is  non-convex due to the constraints
in~\eqref{eqn:up2_quad}, the integrality constraints~\eqref{eqn:up2_bin} and
potentially the policy constraints~\eqref{eqn:up2_policy}. The latter have
been rewritten to depend on $a$ instead of $Y_\text{upg}$ without loss of 
generality since $Y_\text{upg}$ can be computed from $a$ using~\eqref{eqn:yupg}.


\subsection{Voltage magnitude constraints}
A bus voltage magnitude constraint,
\begin{equation} 
  v_{\min,j} \le |\tilde v_j^k| \le v_{\max,j},
\end{equation}
can be rewritten as,
\begin{equation}
  v_{\min,j}^2 \le (v_{r,j}^k)^2 + (v_{q,j}^k)^2 \le v_{\max,j}^2, 
\end{equation}
where the newly introduced variables $v_r^k, v_q^k \in \R^N$ represent the real
and imaginary parts of $\tilde v^k$ (separate vectors for each scenario $k$)
and hence as,
\begin{equation}
  \label{eqn:voltfin}
  v_{\min,j}^2 \le (z^k)^T Q_j z^k \le v_{\max,j}^2,
\end{equation}
where we introduced the shorthand notation $z^k := \bmb (v_r^k)^T & (v_q^k)^T \bme^T$
and where $Q_j$ has entries equal to $1$ in positions $(j,j)$ and $(N+j,N+j)$
and zeros everywhere else. The constraint~\eqref{eqn:voltfin} is now in the
form of~\eqref{eqn:up2_quad}. 

\subsection{Current constraints}
A current constraint,
\begin{equation}
  |Y_{\text{upg},jl}|\cdot|\tilde v^k_j - \tilde v^k_l | \le I_{\max,\text{upg},jl}, \\
\end{equation}
requires slightly more work to reformulate. We first move everything related to
upgrades to the right-hand side and square both sides,
\begin{equation}
  \label{eqn:currupglim}
  |\tilde v_j^k - \tilde v_l^k |^2 \le \frac{I_{\max,\text{upg},jl}^2}
    {|Y_{\text{upg},jl}|^2}.
\end{equation}
We then rewrite the left-hand side as a function of $z$,
\begin{equation}
  \begin{aligned} |\tilde v_j^k - \tilde v_l^k|^2
    &= (v_{r,j}^k - v_{r,l}^k)^2 + (v_{q,j}^k - v_{q,l}^k)^2\\
    &= (z^k)^T Q_{jl} z^k,
  \end{aligned}
\end{equation}
where now $Q_{jl}$ is zero everywhere except entries $1$ at $(j,j), (l,l),
(N+j,N+j), (N+l,N+l)$ and entries $-1$ at $(j,l),(l,j), (N+j,N+l), (N+l,N+j)$.
The right-hand side of~\eqref{eqn:currupglim} depends on the upgrade choices.
Recalling~\eqref{eqn:yupg} and~\eqref{eqn:imaxupg} and the fact that only one
upgrade choice can be made for each line, the fraction
in~\eqref{eqn:currupglim} can be rewritten as an equation that is linear in
$a$,
\begin{equation}
  \begin{aligned} \frac{I_{\max,\text{upg},jl}^2}{|Y_{\text{upg},jl}|^2} &= 
  \frac{I_{\max,jl}^2}{|Y_{jl}|^2}\; + \\
    & \sum_{i \in \mathcal U_{jl}}
    a_i \left[\left(\frac{\Delta I_{jl,i} + I_{\max,jl} }
    {|\Delta (Y_i)_{jl} + Y_{jl}| } 
    \right)^2 - \frac{I_{\max,jl}^2}{|Y_{jl}|^2} \right].
  \end{aligned}
\end{equation}
We can now write the line current constraints as 
\begin{equation}
  \label{eqn:currfin}
  (z^k)^TQ_{jl}z^k + m_{jl}^Ta \le u_{jl},
\end{equation}
where $u_{jl} = I_{\max,jl}^2\,/\,|Y_{jl}|^2$, and 
\begin{equation}
  (m_{jl})_i = -\left[\left(\frac{\Delta I_{jl,i} + I_{\max,jl} }
    {\Delta (Y_i)_{jl} + |Y_{jl}| } 
    \right)^2 - \frac{I_{\max,jl}^2}{|Y_{jl}|^2} \right],
\end{equation}
for $i \in \mathcal U_{jl}$ and $0$ otherwise.
Equation~\eqref{eqn:currfin} now has the same structure as~\eqref{eqn:up2_quad}.

\subsection{Line power constraints}
The concept used to implement the line power constraints is similar to the one
in~\cite{Jabr2013}. We implement Big-M type constraints that select one out of
several equalities depending on which upgrade is selected.  We start from the
power flow equations that each upgrade must satisfy,
\begin{equation}
  \tilde s^k = \diag(\tilde v^k)\overline{Y_{\text{upg}}(a)\tilde v^k}.
\end{equation}
The products of binary and continuous variables in equality constraints
will lead to non-convex relaxations.
In order to avoid such products, we introduce separate case distinctions for
each line as follows:
\begin{equation}
  \label{eqn:linepow_orig}
  \tilde s_{jl}^k = 
  \begin{cases}
    \tilde v_j^k \overline{Y_{\text{upg},jl}(\hat a^1)(\tilde v_l^k-\tilde v_j^k)}, 
      & \text{if } a = \hat a^1 \\
    \tilde v_j^k \overline{Y_{\text{upg},jl}(\hat a^2)(\tilde v_l^k-\tilde v_j^k)}, 
      & \text{if } a = \hat a^2 \\
    \vdots \\
    \tilde v_j^k \overline{Y_{\text{upg},jl}(\hat a^{n_u(jl)})(\tilde v_l^k-\tilde v_j^k)},
      & \text{if } a = \hat a^{n_u(jl)}
  \end{cases}
\end{equation}
where $\tilde s_{jl}^k$ is to be read as ``power flowing into bus $j$, out of
the line $(j,l)$'' and the individual $\hat a^m,\; m \in \{1,\ldots,n_u(jl)\}$, are
the different upgrade possibilities affecting the line $(j,l)$. The~$\tilde s_{jl}^k$
can be expressed using additional variables $f_r^k \in \R^{2L}$ and $f_q^k \in
\R^{2L}$, which represent real and reactive powers flowing into buses from
lines,
\begin{equation}
  \tilde s_{jl}^k = (f_r)_{jl}^k + \sqrt{-1} \cdot(f_q)_{jl}^k.
\end{equation}
In the above, the index $jl$ is used to refer to the entries in $f_r$ and $f_q$
that are assigned to the real and reactive parts of the power flowing into bus
$j$ from the line between buses $j$ and $l$.  
Note that while there can be many binary variables in the complete problem, the
number of variables affecting a particular line is typically small. 
Additional constraints are needed to enforce the power balance for each bus,
\begin{equation}
  \label{eqn:powerbalance}
  \begin{aligned}
    p_{\min,j}^k \le \sum_{l} \real(\tilde s_{jl}^k) \le p_{\max,j}^k, \\
    q_{\min,j}^k \le \sum_{l} \imag(\tilde s_{jl}^k) \le q_{\max,j}^k, 
  \end{aligned}
\end{equation}
where the sum is over all neighboring indices $l$ of $j$.
All constraints for a scenario $k$ of the kind in~\eqref{eqn:powerbalance} are
then collected in the constraints~\eqref{eqn:up2_snappoly}.
As for the actual implementation,
for a line $(j,l)$ and upgrade option $i$ affecting it, we introduce the notation
$Y_{\text{upg},jl}(a_i=1)$ to refer to the admittance of line $jl$ in case $a_i = 1$
(which, by the constraint that just one of the upgrades per line can be
chosen, implies that all other entries of $a$ with indices in~$\mathcal U_{jl}$
are zero).  Using this notation, we would have the constraints
\begin{equation}
  \label{eqn:linepow_lin1}
  \begin{aligned}
    &\left|(f_r^k)_{jl} - \real\left(
    \tilde v_j^k \overline{Y_{\text{upg},jl}(a_i = 1)(\tilde v_l^k-\tilde v_j^k)}
    \right)\right| \le M_{jl}(1-a_i), \\
    &\left|(f_q^k)_{jl} - \imag\left(
    \tilde v_j^k \overline{Y_{\text{upg},jl}(a_i = 1)(\tilde v_l^k-\tilde v_j^k)}
    \right)\right| \le M_{jl}(1-a_i),
  \end{aligned}
\end{equation}
where $M$ is large enough that if $a_i = 0$, the two constraints can never be
active for an otherwise feasible choice of the variables. The choice of $M$ is
important to good relaxation conditioning, which is why it is
discussed in more detail in Section~\ref{ssec:bigm}. In addition to the above,
a ``no upgrade to this line'' case has to be added,
\begin{equation}
  \label{eqn:linepow_lin2}
  \begin{aligned}
    &\left|(f_r^k)_{jl} - \real\left(
    \tilde v_j^k \overline{Y_{jl}(\tilde v_l^k-\tilde v_j^k)}
    \right)\right| \le M_{jl}\sum_{i\in \mathcal U_{jl}} a_i, \\
    &\left|(f_q)_{jl} - \imag\left(
    \tilde v_j^k \overline{Y_{jl}(\tilde v_l^k-\tilde v_j^k)}
    \right)\right| \le M_{jl}\sum_{i\in \mathcal U_{jl}} a_i,
  \end{aligned}
\end{equation}
where the summation over upgrade index $i$ only includes the upgrades that
affect the given line. The quadratic terms in~\eqref{eqn:linepow_lin2} can 
be rewritten as follows:
\begin{equation}
  \label{eqn:linepow_quad}
  \begin{aligned}
  &\real\left( \tilde v_j^k \overline{Y_{jl}(\tilde v_l^k-\tilde v_j^k)} \right) =\\
    &\quad\quad -\real(Y_{jl})(v_r^k)_j^2 + \real(Y_{jl})(v_r^k)_j(v_r^k)_l \\
    &\quad\quad - \imag(Y_{jl}) (v_r^k)_j(v_q^k)_l + \real(Y_{jl})(v_q^k)_j(v_q^k)_l \\
    &\quad\quad - \real(Y_{jl}) (v_q^k)_j^2 + \imag(Y_{jl})(v_q^k)_j(v_r^k)_l, \\
  &\imag\left( \tilde v_j^k \overline{Y_{jl}(\tilde v_l^k-\tilde v_j^k)} \right) = \\
    &\quad\quad \real(Y_{jl})(v_q^k)_j(v_r^k)_l - \imag(Y_{jl})(v_q^k)_j(v_q^k)_l \\
    &\quad\quad + \imag(Y_{jl})(v_q^k)_j^2 - \real(Y_{jl})(v_r^k)_j(v_q^k)_l \\
    &\quad\quad - \imag(Y_{jl})(v_r^k)_j(v_r^k)_l + \imag(Y_{jl})(v_r^k)_j^2.
  \end{aligned} 
\end{equation}
A similar procedure can be used for~\eqref{eqn:linepow_lin1}, simply by replacing
$Y_{jl}$ with $Y_{\text{upg},jl}(a_i = 1)$.
This means the line power constraints are now quadratic in~$z^k$ and linear 
in~$y^k := \bmb (f_r^k)^T & (f_q^k)^T \bme^T$ and $a$, as was required (the
absolute value operator can simply be replaced with two linear constraints).
The data for the $Q$ terms is given in~\eqref{eqn:linepow_quad}. The data for
the $q$ terms is determined by how the variables are ordered in $f_r^k$ and
$f_q^k$ and by~\eqref{eqn:linepow_lin1} and~\eqref{eqn:linepow_lin2}.  The data for
the $m$ terms is also determined
by~\eqref{eqn:linepow_lin1} and~\eqref{eqn:linepow_lin2}.  Each set of constraints
in~\eqref{eqn:linepow_lin1} or~\eqref{eqn:linepow_lin2} translates into 4
constraints of the form~\eqref{eqn:up2_quad} due to the absolute value involved. 

\subsection{Computation of Big-M terms for line powers}
\label{ssec:bigm}
Problem~\eqref{eqn:upgradeproblem2} will later be solved in a Branch-and-Bound
setting by relaxing the integrality constraints on $a_i$ to be $a_i \in [0,1]$.
The constant $M$ should be chosen large enough
for~\eqref{eqn:linepow_lin1}-\eqref{eqn:linepow_lin2} to
implement~\eqref{eqn:linepow_orig}, but should
also be as small as possible to avoid numerical issues. For this reason,
instead of a single $M$ for all constraints, a separate constant $M_{jl}$ is
found for each line from bus $j$ to bus $l$.
The following Lemma gives a lower bound on~$M_{jl}$.
\begin{lemma}
In order for the case distinction~\eqref{eqn:linepow_orig} to be equivalent to
the intersection of the constraints in~\eqref{eqn:linepow_lin1}
and~\eqref{eqn:linepow_lin2}, it has to hold that
\begin{equation}
  \label{eqn:bigmlb1}
  \begin{aligned}
    &M_{jl} \ge v_{\max,j}
      \max \left| \tilde v_l^k - \tilde v_j^k \right| \cdot \\
    &\quad \max_{i_1,i_2 \in \mathcal U_{jl}} 
      \big|Y_{\text{upg},jl}(a_{i_1}=1)-Y_{\text{upg},jl}(a_{i_2}=1)  \big|, \\
  \end{aligned}
\end{equation}
where the term $\max |\tilde v_l^k - \tilde v_j^k|$ refers to the maximum
absolute value of the difference between $\tilde v_l^k$ and $\tilde v_j^k$ that
can occur in scenario $k$.
\end{lemma}
\begin{proof}
For any integral choice of $a$ satisfying~\eqref{eqn:upgcon}, one of the cases
in~\eqref{eqn:linepow_orig} is selected by means of the right hand side of one
of the corresponding constraints
in~\eqref{eqn:linepow_lin1}--\eqref{eqn:linepow_lin2} becoming $0$. This is
referred to hereafter as this constraint being \emph{active}.  What is to
be shown is that if $M_{jl}$ is chosen to satisfy~\eqref{eqn:bigmlb1}, none of
the other constraints in~\eqref{eqn:linepow_lin1}--\eqref{eqn:linepow_lin2} can
be violated for any admissible choice of $\tilde v^k$. For this, the largest
absolute value that any of the non-active constraints can attain is to be
found and then $M_{jl}$ has to be picked larger than that.  We first define a
  shorter version of the notation introduced after~\eqref{eqn:powerbalance}, 
\[ y_{jl}(i) := Y_{\text{upg},jl}(a_i = 1), \]
as well as a shorthand notation for the difference between voltages of two
buses $j$ and $l$,
$\Delta \tilde v_{jl}^k := \tilde v_j^k - \tilde v_l^k$.
We can then write the largest possible absolute value as
\begin{equation}
  \begin{aligned}
    &\max_{i_1,i_2} \max_{\tilde v^k} \Big| 
      \real\left( \tilde v_j^k \overline{y_{jl}(i_1)\Delta \tilde v_{lj}^k} \right) 
      -\real\left( \tilde v_j^k \overline{y_{jl}(i_2)\Delta \tilde v_{lj}^k}\right)\Big|\\
      &\le \max_{i_1,i_2}\max_{\tilde v^k} \Big| 
      \tilde v_j^k \overline{\left(y_{jl}(i_1)-y_{jl}(i_2)\right)
      \Delta \tilde v_{lj}^k} \Big|  \\
    &= \max_{\tilde v^k} \left|\tilde v_j^k \right|
      \max_{\tilde v^k} \left|\Delta \tilde v_{lj}^k\right|
      \max_{i_1,i_2} \Big| \left(y_{jl}(i_1)-y_{jl}(i_2)\right)\Big|.
  \end{aligned}
\end{equation}
The first term in the last row above is just $v_{\max,j}$ from the problem data.
The last line of the above therefore equals~\eqref{eqn:bigmlb1}, which
completes the proof. 
\end{proof}

The last term in~\eqref{eqn:bigmlb1} can easily be evaluated exactly by
enumeration of the possible pairs of different line parameters. 
A bound on $\Delta \tilde v_{lj}^k$ is supplied by the current limits
through~\eqref{eqn:currupglim}. 
In order to get a bound that holds for all
possibilities, a maximization over $a$ is performed.
\begin{equation} 
  \label{eqn:bettercurrbnd}
  \Delta \tilde v_{lj}^k = |\tilde v_l^k - \tilde v_j^k |^2 \le \max_{a} 
    \frac{I_{\max,\text{upg},lj}^2(a)} {|Y_{\text{upg},lj}(a)|^2}.
\end{equation}

\section{Algorithm}
\label{sec:algorithm}
In order to solve~\eqref{eqn:upgradeproblem2}, a Branch-and-Bound procedure
augmented with constraint generation is applied to its mixed-integer
semidefinite relaxation. This relaxation has the form
\begin{subequations}
  \makeatletter
  \def\@currentlabel{R}
  \makeatother
  \label{eqn:upgradeproblemrel}
  \renewcommand{\theequation}{R.\arabic{equation}}
  \vspace{-0.6cm}
  \begin{align}
    \text{Problem R:} \nonumber \\ 
    \minim_{\substack{a,Z^k,y^k,\\ k \in \{1, \ldots, K\}}} \;\;& 1^Ta \\
         \st \;\;& a \in \{0,1\}^{n_u}\label{eqn:uprel_bin}, \\ 
             & Aa \le b, \label{eqn:uprel_poly}\\
             & Cy^k \le d^k, \label{eqn:uprel_snappoly}\\
             & \alpha_h \le \mathrm{tr}(Q_hZ^k) + q_h^Ty^k + m_h^Ta \le \beta_h,
               \label{eqn:uprel_quad} \\
             & Z^k \succeq 0, \\
             & \forall h \in \{1,\ldots, H\}, \forall k \in \{1,\ldots, K\}.
               \nonumber
  \end{align}
\end{subequations}
Problem~\eqref{eqn:upgradeproblemrel} is a relaxation
of~\eqref{eqn:upgradeproblem2}: The policy constraints are removed and
the non-convex quadratic constraints are replaced by a semidefinite
relaxation.
This algorithm was presented in~\cite{merkli2018globally} along
with a proof that it solves~\eqref{eqn:upgradeproblem2}. It
will be extended here to work in concert with the Benders decomposition. 

\subsection{Benders decomposition}
\label{ssec:benders}
Note that the convex problems formed by relaxing the integrality constraints
in~\eqref{eqn:upgradeproblemrel} require one semidefinite variable~$Z^k$ as
well as one vector of auxiliary variables~$y^k$ for each scenario~$k$.
The full set of power flow constraints and the operational constraints also
need to be added separately for each of these scenarios. For $K \gg 1$,
problem~\eqref{eqn:upgradeproblemrel} becomes difficult to handle in
a centralized manner. 

The generalized Benders decomposition~\cite{Benders1962,Geoffrion1972} provides
an approach for dealing with such structured problems efficiently. For a problem 
of the form
\begin{equation}
  \label{eqn:bendersproto}
  \begin{aligned}
    \minim_{\xi,\zeta}\;\; &  F(\zeta) \\
    \st\;\; & G(\xi,\zeta) \le 0,\; \xi \in \mathcal X, \zeta \in \mathcal Z, 
  \end{aligned}
\end{equation}
with $F,G,\mathcal X,\mathcal Z$ convex, the method provides an algorithm that
alternately fixes one of~$\zeta$ and~$\xi$ while solving the problem for the
other.  The minimization with respect to~$\zeta$ is referred to as the ``master
problem'', defined as follows:
\begin{subequations}
  \makeatletter
  \def\@currentlabel{M}
  \makeatother
  \label{eqn:bendersmaster}
  \renewcommand{\theequation}{M.\arabic{equation}}
  \begin{align}
    \text{Problem M:} \nonumber \\ 
    \minim_{\zeta}\;\; & F(\zeta) \\
    \st\;\; & \min_{\xi\in \mathcal X}\left\{ \lambda^T 
        G(\xi,\zeta)\right\} \le 0, \label{eqn:bendersmaster_con}\quad
      \forall \lambda \in \Lambda, \\
      & \zeta \in \mathcal Z.
  \end{align}
\end{subequations}
where $\Lambda$ is a set of vectors which will be iteratively built in the
algorithm.  Conversely, the minimization for $\xi$ is referred to as the
``subproblem'' and defined as
\begin{subequations}
  \makeatletter
  \def\@currentlabel{S}
  \makeatother
  \label{eqn:benderssub}
  \renewcommand{\theequation}{S.\arabic{equation}}
  \begin{align}
    \text{Problem S:} \nonumber \\ 
    \minim_{\xi,\gamma \in \R}\;\; & \gamma \\
    \st\;\; & G(\xi,\zeta) \le 1\gamma, \label{eqn:benderssub_con} \\
      & \xi \in \mathcal X.
  \end{align}
\end{subequations}
The full Benders decomposition algorithm is now stated in
Figure~\ref{alg:benders}, and the interested reader is referred
to~\cite{Benders1962,Geoffrion1972} for more discussion. Each iteration of the
algorithm either finds the optimal solution to~\eqref{eqn:bendersproto} or
improves the lower bound on the optimal cost. The Benders iteration can hence
be stopped at any point and the lower bound obtained at that point is valid.
This makes the algorithm particularly suitable for solving relaxations in a
Branch-and-Bound procedure.
In our application, we partition the variables
of Problem~\eqref{eqn:upgradeproblemrel} as 
\begin{equation}
    \xi^k := (Z^k,y^k),  \quad \zeta := a.
\end{equation}
We refer to the collection of $\xi^k$ for~$k \in \{1,\ldots,K\}$~as $\xi$.
We can then define the sets $\mathcal X$ and $\mathcal Z$:
\begin{equation}
  \begin{aligned}
    \mathcal X^k &:= \left\{ 
      (Z^k,y^k) \;\middle|\; 
      \begin{aligned} 
        Cy^k &\le d^k\\
        Z^k &\succeq 0
      \end{aligned} 
    \right\}, \quad \forall k \in \{1,\ldots,K\}, \\
    \mathcal X &:= \bigcap_{k=1}^K \mathcal X^k,\\
  \end{aligned}
\end{equation}
and
\begin{equation}
  \begin{aligned}
    \mathcal Z &:= \left\{ 
      a \;\middle| \;
      \begin{aligned} 
        &a \in [0,1]^{n_u}\\
        &Aa \le b
      \end{aligned}
    \right\}.
  \end{aligned}
\end{equation}
Additionally, those constraints in~\eqref{eqn:uprel_quad} that correspond 
to the voltage magnitude constraints~\eqref{eqn:opconst1} can also be added to 
$\mathcal Z$ since they do not depend on $a$. The rest of the constraints 
in~\eqref{eqn:uprel_quad} are now stacked to make $G^k(\xi,\zeta)$: 
\begin{equation}
  G^k(\xi^k,\zeta) := \begin{bmatrix}
    -\mathrm{tr}(Q_{1}Z^k) - q_{1}^Ty^k - m_1^Ta + \alpha_1 \\[0.1cm]
    \mathrm{tr}(Q_{1}Z^k) + q_{1}^Ty^k + m_1^Ta - \beta_1 \\[0.1cm]
    -\mathrm{tr}(Q_{2}Z^k) - q_{2}^Ty^k - m_2^Ta + \alpha_2 \\[0.1cm]
    \mathrm{tr}(Q_{2}Z^k) + q_{2}^Ty^k + m_2^Ta - \beta_2 \\
    \vdots
  \end{bmatrix},
\end{equation}
and the $G^k(\xi^k,\zeta)$ are in turn stacked to make $G(\xi,\zeta)$.
We note here that the constraint data $Q_h, q_h, m_h, \alpha_h$ and $\beta_h$ are
the same for all scenarios. Because the entries of $G(\xi,\zeta)$ are linear in
$a$, the maximization for $\xi$ in~\eqref{eqn:bendersmaster_con} can be
performed independently of $\zeta$ by solving a semidefinite problem. The
resulting constraints are linear in $\zeta$ and independent of $\xi$,
making~\eqref{eqn:bendersmaster} a linear problem in $\zeta$. 
Subproblem~\eqref{eqn:benderssub} becomes a semidefinite problem due to
$\mathcal X$. 

\subsubsection{Decomposition of Benders subproblem} 
The main computational burden lies in solving the
subproblem~\eqref{eqn:benderssub} as well as the parametric maximizations
in~\eqref{eqn:bendersmaster_con}. We will now discuss how these two problems can
be solved in a separable manner.  Consider the~$K$ individual problems,
\begin{subequations}
  \makeatletter
  \def\@currentlabel{S$^k$}
  \makeatother
  \label{eqn:benderssubK}
  \renewcommand{\theequation}{S$^k$.\arabic{equation}}
  \begin{align}
    \text{Problem S$^k$:} \nonumber \\ 
    \minim_{\xi^k,\gamma^k \in \R}\;\; & \gamma^k \\
    \st\;\; & G^k(\xi^k,\zeta) \le 1\gamma^k, \label{eqn:benderssub_conK} \\
      & \xi^k \in \mathcal X^k.
  \end{align}
\end{subequations}
The optimal results $(\gamma^k)^*, (\xi^k)^*$ of these problems can be used to
construct an optimal solution for problem~\eqref{eqn:benderssub}. It holds that
\begin{equation}
  \gamma^* = \max_k\; (\gamma^k)^*
\end{equation}
where $\gamma^*$ denotes the optimal $\gamma$ for~\eqref{eqn:benderssub}.
The~$(\xi^k)^*$ are feasible for~\eqref{eqn:benderssub} with the above choice
of~$\gamma^*$ as well.  The dual multipliers, however, are different. The
complementarity conditions for constraint~\eqref{eqn:benderssub_con} can be
written as 
\begin{equation}
  \label{eqn:kktbenderssub}
  \lambda_k^T\Big(G^k(\xi^k,\zeta)-1\gamma\Big) = 0, \quad k \in \{1,\ldots, K\}
\end{equation}
due to all the involved terms being non-negative. Hence for a given choice of
optimal variables, it will hold that~$\lambda_k = 0$ for all constraints for
which $G^k((\xi^k)^*,\zeta) < 1\gamma^*$.  Considering that $(\gamma^k)^* \ge \gamma^*$, 
let~$\mathcal K$ be defined as 
\begin{equation} 
  \mathcal K := \left\{ k \in 1,\ldots,K \; \middle| \; G^k((\xi^k)^*,\zeta) = \gamma^* \right\}.
\end{equation}
For the optimal multipliers $\lambda^k$ of~\eqref{eqn:benderssub}, it holds that
\begin{equation}
  \lambda^k = 0 \quad \forall k \not\in \mathcal K
\end{equation}
due to~\eqref{eqn:kktbenderssub}. In order to obtain the $\lambda_k$ for the $k
\in \mathcal K$, a reduced version of~\eqref{eqn:benderssub} can be solved
where $G$ is constructed from only the scenarios in $\mathcal K$. Note that in
practice, $|\mathcal K| \ll |K|$, in fact it is exceedingly unlikely that
$|\mathcal K| > 1$ occurs. If~$|\mathcal K| = 1$, the individual $\lambda_k$
for the single $k$ in $\mathcal K$ can be used as is,
with all the others set to $0$.

\subsubsection{Decomposition of the parametric cut problem}
The constraints~\eqref{eqn:bendersmaster_con} are parametrically solved in
$\xi$ once for each $\lambda$ added to $\Lambda$. Due to the structure
of $G$, the parametric problem can be written as follows:
\begin{equation}
  \begin{aligned}
    \min_{\xi\in \mathcal X}\left\{ \lambda^T G(\xi,\zeta)\right\} =&  
      \min_{\xi^k\in \mathcal X^k}\left\{ \sum_{k=1}^K \lambda_k^T G^k(\xi^k,\zeta)\right\} \\
      =& \sum_{k=1}^K \min_{\xi^k\in \mathcal X^k}\left\{\lambda_k^T G^k(\xi^k,\zeta)\right\},
  \end{aligned}
\end{equation}
which means the parametric cut problems can be solved separately for each $k$
and the results summed up.



\begin{figure}
  \vspace{0.2cm}
  \begin{algorithmic}[1]
    \small
    \State Initialize $\Lambda = \emptyset$
    \While{No feasible solution $\zeta$ found}
      \State Solve~\eqref{eqn:bendersmaster} for $\zeta$ (solution: $\zeta^\star$)
      \State Solve~\eqref{eqn:benderssub} for $\xi$ with $\zeta$ fixed to $\zeta^\star$
       (solution: $\xi^\star,\gamma^\star$)
      \If{If $\gamma^\star \le 0$}
        \State Stop, $\xi^\star,\zeta^\star$ are optimal for~\eqref{eqn:bendersproto}
      \Else
        \State Let $\lambda^\star$ be the dual multipliers of constraints~\eqref{eqn:benderssub_con}
        \State Add $\lambda^\star$ to $\Lambda$
      \EndIf
    \EndWhile
  \end{algorithmic}
  \caption{Benders decomposition algorithm for problems of the
    form~\eqref{eqn:bendersproto}. Intuitively speaking, the algorithm picks
    the best $\zeta$ possible without considering $\xi$
    using~\eqref{eqn:bendersmaster}, then attempts to find a $\xi$ that is
    still feasible using~\eqref{eqn:benderssub}. If that fails, the information
    from~\eqref{eqn:benderssub} is used to restrict the search
    space in~$\zeta$.} 
  \label{alg:benders}
\end{figure}

\subsection{Combined Branch-and-Bound and Benders algorithm}
\begin{figure}
  \vspace{0.2cm}
  \begin{algorithmic}[1]
    \small
    \State Set $U = \infty, L = -\infty$, 
      Tree: Root vertex $\mathcal I_0 = \mathcal I_1 = \emptyset$ 
    \label{algstep:termination}
    \While{$U - L > \varepsilon$}
      \State Pick an unprocessed vertex $\mathcal N$ with index
      sets $\mathcal I_0^{\mathcal N}$, $\mathcal I_1^{\mathcal N}$
      \State Perform at most $B$ benders iterations
        on~\eqref{eqn:upgradeproblemrel} with~\eqref{eqn:uprel_bin} replaced by 
      \[ a_i \in \begin{cases} \{0\}, & \text{if } i \in \mathcal I_0^{\mathcal N}, \\
        \{1\}, & \text{if } i \in \mathcal I_1^{\mathcal N}, \\
          [0,1] & \text{otherwise.} \end{cases} 
      \] and the Benders cuts from all ancestors added\label{alg:solveref}
      \If{Problem in step~4 was not infeasible}
        \State Store the obtained Benders cuts in this node
        \State Let $(a^{\mathcal N}, Z^{\mathcal N}, y^{\mathcal N})$ refer to the Benders result 
        \State Let $L^{\mathcal N} = 1^Ta^\mathcal N$
        \If{Feasible, $L^{\mathcal N} < U$ and $a^{\mathcal N} \in \{0,1\}^{n_u}$} 
          \vspace{0.1cm}
          \State Evaluate policy $g(a^\mathcal N, \mathrm{limits}^k), \; \forall k$ 
            \label{alg:policycheck}
          \vspace{0.1cm}
          \If{Feasible for all $k$}
            \vspace{0.1cm}
            \State Update $U = 1^Ta^{\mathcal N}$ \label{alg:updateu}
            \vspace{0.1cm}
          \Else
            \vspace{0.1cm}
            \State Add cut $\|a - a^{\mathcal N}\|_1 \ge 1$
            \State Go back to the solve step (line \ref{alg:solveref})
            \vspace{0.1cm}
          \EndIf \label{alg:policyend}
          \vspace{0.1cm}
        \ElsIf{$L^{\mathcal N} < U$ but $a^{\mathcal N} \not\in \{0,1\}^{n_u}$}
          \vspace{0.1cm}
          \State Select index $\ell$, $\ell \not\in \mathcal I_0^{\mathcal N} 
            \cup \mathcal I_1^{\mathcal N}$
          \vspace{0.1cm}
          \State Add a vertex with $\mathcal I_0 = \mathcal I_0^{\mathcal N} 
            \cup \{\ell\}, \mathcal I_1 =\mathcal  I_1^{\mathcal N}$
          \vspace{0.1cm}
          \State Add a vertex with $\mathcal I_0 = \mathcal I_0^{\mathcal N} ,
            \mathcal I_1 = \mathcal I_1^{\mathcal N}\cup \{\ell\}$
          \vspace{0.1cm}
        \EndIf 
      \Else\vspace{0.1cm}
        \State Set lower bound for this subtree to $\infty$
        \vspace{0.1cm}
      \EndIf \vspace{0.1cm}
      \State Update $L = \min \left\{ L^{\mathcal N}\; |\; \mathcal N \in \text{tree} \right\}$ 
      \State Run heuristic methods
      \If{Heuristic methods yield feasible solution}
        \State Update $U$ to reflect best solution found
      \EndIf
      \vspace{0.1cm}
    \EndWhile
  \end{algorithmic}
  \caption{Branch-and-Bound algorithm with policy constraint generation. The algorithm
    traverses the binary tree based on fixing entries $a_i$ to either $1$ or $0$, relaxing
    the non-fixed entries to $[0,1]$. The first difference between regular
    Branch-and-Bound and the algorithm here is the policy evaluation and
    constraint addition. See~\cite{merkli2018globally} for further details.
    The second difference is the application of incomplete Benders iterations
    for the relaxations, which is discussed in the text.} 
  \label{alg:bnbcuts}
\end{figure}
The modified Branch-and-Bound algorithm is shown in Figure~\ref{alg:bnbcuts},
extended to include the changes introduced by the application of the Benders
decomposition for the relaxations. The major differences to the version
presented in~\cite{merkli2018globally} are as follows:
\begin{enumerate}[1)]
  \item The relaxations are solved with the Benders decomposition. Only a
    limited number, $B \in \N$, of Benders iterations are performed on each
    relaxation.
  \item The algorithm reuses cuts obtained from the partial solves for
    tightened versions of the relaxation from which they were obtained, since
    they remain valid.
\end{enumerate}
Heuristic methods include approaches such as the ones outlined in the
introduction. Well-performing heuristics can lead to good or even optimal
solutions being found very rapidly, even though the solution process will
typically take a much longer time to certify that an optimal solution is
actually optimal.

The following Lemmas provide 
an explanation as to why this combination of algorithms is computationally 
correct in the sense that no feasible points are artificially discarded.
\begin{lemma}
  \label{lem:lowerbound}
  Despite not solving the relaxed problems fully, the algorithm never cuts
  branches of the tree that would not have been cut if the relaxations were
  fully solved.
\end{lemma}
\begin{proof}
  Each Benders iteration yields a valid lower bound on the objective,
  $L_{\text{Benders}} \le L_\text{true}$.
  This means that for the result of each iteration, 
  \[ L_{\text{Benders}} > U \implies L_\text{true} > U, \]
  and hence the branch cutting criterion can never discard a branch that should
  have been explored.
\end{proof}
\begin{lemma} 
  \label{lem:cuts}
  Cuts obtained in each Benders solve are valid for all
  descendants of that node in the Branch-and-Bound tree.
\end{lemma}
\begin{proof}
  The feasible set of each node is a restriction of that of all its parents.
  All constraints that were valid for the parent feasible set are also valid
  for the more restricted feasible set of the node itself.
\end{proof}

Lemma~\ref{lem:cuts} does not make a statement about whether the lower bound of
the relaxations eventually reaches the true lower bound. While this is not
explored further in this work, the numerical experiments suggest that the
obtained lower bounds approximate the true values well enough for the algorithm
to be effective.

\section{Numerical experiments}
\label{sec:numexp}
For the software implementation of the experiments in this section, the Julia
language~\cite{Bezanson2017} was used in conjunction with the JuMP modeling
package~\cite{DunningHuchetteLubin2017}.  Semidefinite relaxations were solved
with MOSEK, smooth nonlinear problems arising in the operating policy with
IPOPT~\cite{Waechter2005} and linear problems with Gurobi. The experiments were 
run on a machine with 64~GB of system memory and an AMD Ryzen~9 3950X 
(3.5~GHz). All experiments were performed in single-threaded mode.
The operating system used was Debian Linux~10.

\subsection{Experiment considerations}
Both the centralized problem as well as the optimization problems in the Benders-based
approach were solved using a single thread. This was done to make
comparisons easier, such as comparing a single thread centralized solution to the Benders
approach with 8 parallel workers. In practice, most numerical solver codes such as MOSEK
can make use of multiple threads on the same machine and such parallelism can trivially be
combined with the parallelism provided by the method presented here.

As mentioned in Lemma~\ref{lem:lowerbound}, the number of Benders iterations performed for
a single relaxation can be chosen freely. Setting the iteration number high enough 
to ensure that all relaxations are solved to a high accuracy was deemed impractical since
a large number of iterations would be required. On the other hand, setting a low
maximum number of iterations leads to the lower bounds on the relaxation objective value
being lower than the true relaxation objective. This in turn leads to an increase in the 
number of Branch-and-Bound nodes visited, because the tree can be cut less often.
Experiments have shown that as few as between 3 and 10 Benders iterations per 
relaxation provide a good tradeoff between the time taken for each relaxation and the 
increase in the number of nodes visited by the Branch-and-Bound procedure due to the 
lower bounds being weaker.

One can hence say that the Branch-and-Benders approach makes sense in cases where the
computation speed gained by increased parallelism is larger than the speed lost due to
an increase in the number of nodes visited. While this statement is straightforward to
understand and even quantify on a single problem instance, it is relatively challenging 
to generalize to a set of rules for when the method presented here makes sense.
The approach taken here is to present an example that was small enough to run in a
reasonable time frame on a single machine over a range of scenario numbers, discuss the 
results of this particular experiment and make some cautious extrapolation attempts 
based on the data.

\subsection{Test problem instance}
\begin{figure}
  \centering
  \includegraphics[width=0.7\columnwidth]{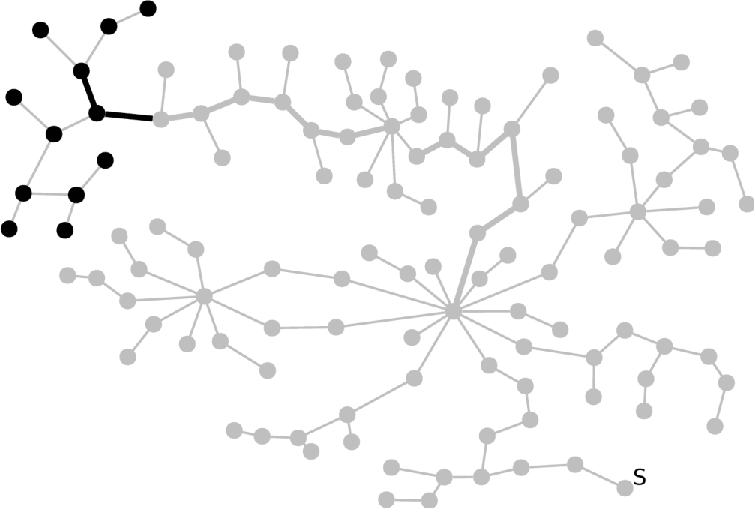}
  \caption{Visualization of part of the Zurich distribution grid and violations
  encountered with a simulated load scenario. Dark vertices represent buses
  with voltage magnitude violations, dark lines represent current limit
  violations. The slack bus is denoted by ``S''. The thicker lines represent
  lines that are considered for upgrades in the study.}
  \label{fig:EKZvio}
\end{figure}
A realistic industrial case study was performed with part of the Zurich
distribution grid. Actual load data was augmented with some simulated PV
in-feeds. Additional studies on this data are given in~\cite{merkli2018globally}.  Load
data as well as the system topology data are the same as in~\cite{merkli2018globally}, 
but only the strongest of the upgrade possibilities for a subset of lines is considered 
here to shorten computation times. A visual representation of the violations encountered 
for the system and load data is shown in Fig.~\ref{fig:EKZvio}. Different cardinality 
sets of randomly generated scenarios were used in the following experiments. The scenarios
were created by perturbing a base scenario obtained from real load data. In the experiment
presented here, the policy was simply solving an AC economic dispatch to local optimality
using IPOPT.

A greedy heuristic based on adding one upgrade at a time was used to look for feasible
solutions. This heuristic was run on every feasible relaxation encountered and was
computationally rather expensive, as for each upgrade addition the operating policy
is evaluated for each scenario separately. The greedy heuristic was also trivially
parallelizable, however. 

\subsection{Runtime comparison}
\begin{table}[t]
  \centering
  \vspace{-0.1cm}
  \renewcommand{\arraystretch}{1.3}
  \caption{Comparison of centralized and Benders solve times}
  \begin{tabular}{l  r r r}
  \textbf{1 scenario} & \textbf{centralized} & \textbf{Benders 3/1} & \textbf{Benders 3/K}  \\
  \hline
Final objective:     &          3 &                       3 &          3 \\
Nodes visited:       &          6 &                      14 &         14 \\
Policy cuts:         &          1 &                       1 &          1 \\
Policy checks [s]:   &        1.8 &                     1.8 &        1.8 \\
Heuristics [s]:      &        8.5 &                    19.5 &       19.5 \\
Relaxations [s]:     &       18.8 &       56.8 +       29.0 &       56.8 +       29.0 \\
Total time [s]:      &       29.1 &                   107.1 &      107.1 \\
  \end{tabular}\\[0.2cm]
  \begin{tabular}{l  r r r}
  \textbf{20 scenarios} & \textbf{centralized} & \textbf{Benders 3/1} & \textbf{Benders 3/K}  \\
  \hline
Final objective:     &          4 &                       4 &          4 \\
Nodes visited:       &         42 &                      86 &         86 \\
Policy cuts:         &         36 &                      14 &         14 \\
Policy checks [s]:   &      693.9 &                   277.0 &       15.2 \\
Heuristics [s]:      &     1460.4 &                  3333.5 &      183.3 \\
Relaxations [s]:     &     2043.1 &     6974.6 +      145.5 &      383.6 +      145.5 \\
Total time [s]:      &     4197.5 &                 10730.6 &      727.7 \\
  \end{tabular}
  \label{tbl:comparison_approaches}
  \vspace{-0.4cm}
\end{table}

The runtimes for the centralized and Branch-and-Benders approaches are shown
in Table~\ref{tbl:comparison_approaches}. The column named ``Benders 3/1'' refers to 
a configuration where at most 3 Benders iterations were performed on each relaxation
solve and everything runs in a single-threaded configuration. Similarly, ``Benders 3/K'' 
refers to a hypothetical implementation where $K$ workers are used in parallel. The
relaxation times are separated into two numbers in the Benders columns, the first being
the subproblem time and the second being the sum of the master and cut problem times. For
the computation of the runtimes in the parallel scenario, the subproblem time was divided
by the number of agents $K$ and 10\% was added to the result to account for communication
overhead that a real implementation would have. The cut and master problem times were
taken as they are, because no parallelisation is suggested in the method here for them.
Similarly, the time for the heuristics and policy checks was divided by the number of
agents and 10\% was added to the resulting time for coordination.

Investigating the data, one can see that for a single scenario the Benders-based approach
makes little sense --- the overhead that stems from having to solve multiple Benders
iterations combined with more nodes being visited leads to a much worse overall runtime.
For 20 scenarios, the situation is different. While taking the Branch-and-Benders approach
with a single thread still takes more than twice the time as solving the centralized
problem, the Benders approach scales very well with the number of parallel processing
agents. For the best case of having as many processing agents as there are scenarios, the
Benders approach winds up using less than 18\% of the runtime of the centralized approach,
despite visiting more than twice as many nodes. 

This result is of course optimistic in some aspects: The centralized approach is expected
to also scale well if the numerical solver used is given more threads to work with, and
the parallelization of both heuristics and policy checks can be done independently of how
the relaxations are solved. Regardless, a problem of real scale may have thousands of
scenarios, growing the gap between a centralized solution with a many-core machine and the
Benders approach with dozens or even hundreds of machines at its disposal. Note that this
is not because the Benders approach is inherently more computationally efficient, but
because it allows the computation to be distributed across a large set of computers.
As mentioned earlier, additional forms of parallelization such as solving individual 
subproblems with multiple threads or examining multiple Branch-and-Bound tree nodes
concurrently can be combined with what is presented here without any issues.

\subsection{Discussion of runtime behavior}
\begin{figure}
  \centering
  \includegraphics[width=0.7\columnwidth]{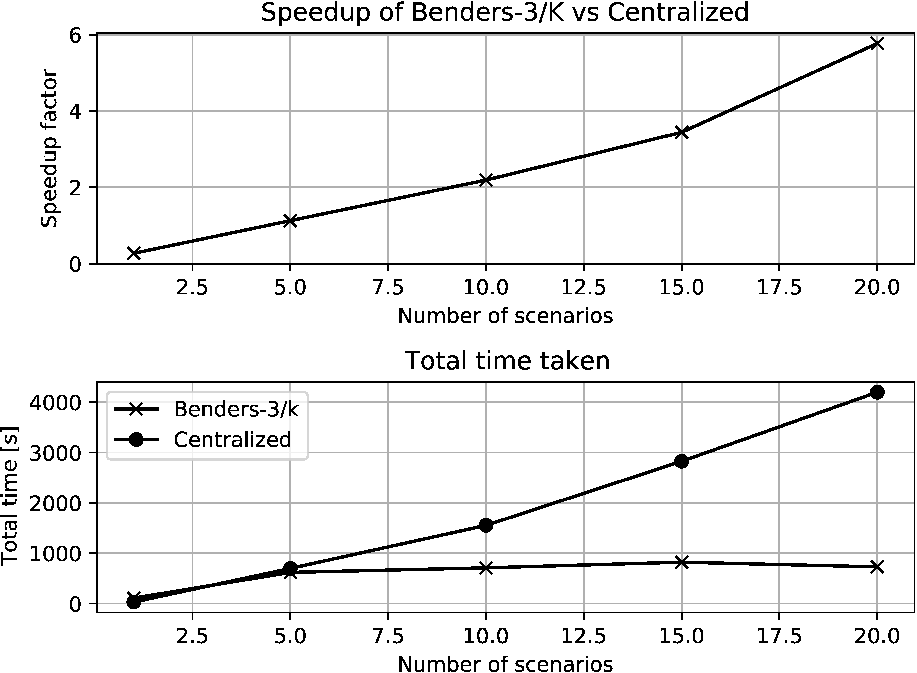}
  \caption{Comparison of complete solve runtimes of centralized and fully parallel Benders approach for a
  series of scenario numbers.}
  \label{fig:scaling_with_scenarios}
\end{figure}
In order to aid potential users of the method decide whether it makes sense for their
problem scale, this section will discuss the expected behavior of runtimes based
on small-scale experiments with 1, 5, 10, 15 and 20~scenarios. In Figure~\ref{fig:scaling_with_scenarios},
graphs are shown comparing the first and third solution methods listed in Table~\ref{tbl:comparison_approaches} for
these numbers of scenarios. 

In the upper plot, it can be seen that the computational advantage steadily increases with the number of 
scenarios. This behavior is expected, because for the best case of as many agents as there are scenarios,
the time per Benders iteration is roughly constant, whereas the centralized solution time grows at least
linearly with the number of scenarios. In the lower plot, total runtimes are shown for comparison, though
care has to be taken when showing absolute times for different scenario numbers in the same plot. This is 
because the Branch-and-Bound trees, number of heuristic evaluations and number of nodes visited vary across
the different scenario counts.

These experiments showed quite regular behavior, because they all required only 3 to 4 upgrades to the system
and hence required only a slowly growing number of nodes to be visited as more scenarios were added. If more
scenarios are added or a different power system is optimized over, the number of upgrades required can become
higher and with that the runtime of both the centralized and parallel approaches can grow substantially. 
However, as long as the number of nodes required for Benders does not grow strongly compared to the number 
of nodes required for the centralized approach, the advantage gained due to the Benders-based parallelization
is expected to persist and even grow with an increasing number of scenarios. 

\begin{figure}
  \centering
  \includegraphics[width=0.7\columnwidth]{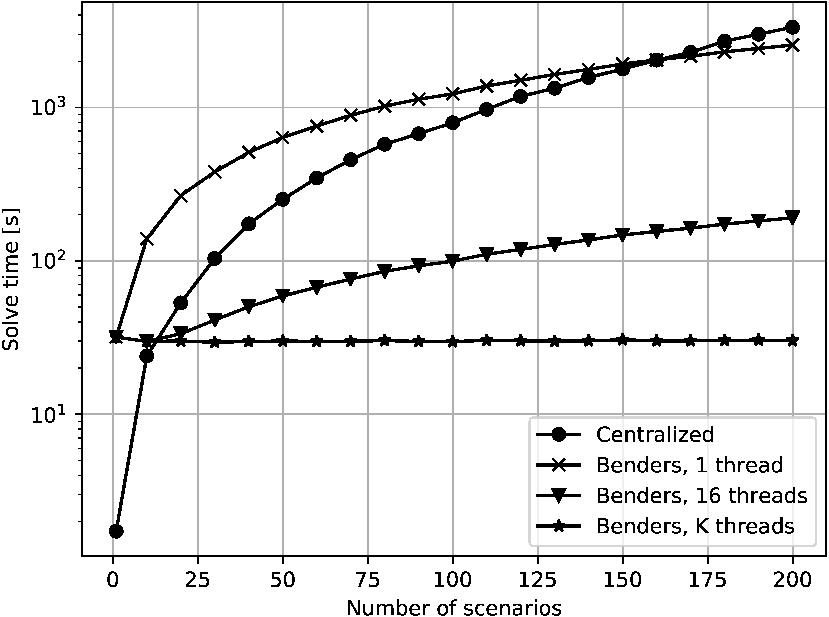}
  \caption{Comparison of runtimes of a single relaxation using the centralized, 
  Benders with 16 parallel computing agents and fully parallel ($K$ agents) 
  Benders methods}
  \label{fig:single_scaling_with_scenarios}
\end{figure}
In Figure~\ref{fig:single_scaling_with_scenarios}, timings are shown for a single relaxation instead, but up
to 200 scenarios. The Benders times were scaled to 10 iterations for the experiment. For a problem with 200
scenarios, 10 iterations of Benders implementation with 16 parallel agents would outperform the centralized
approach by more than an order of magnitude. For a fully parallel implementation, the performance gap would
even grow to two orders of magnitude.

\section{Conclusion}
This work presents an algorithmic framework for effectively solving power
system line upgrade problems at a scale applicable to many city distribution
networks. The method deterministically finds globally optimal
solutions provided they exist, and certificates that they do not otherwise. The
use of the Benders decomposition for solving relaxations renders the core
computational burden of the algorithm parallelizable to a large extent.

While more research is needed regarding the impact of weaker lower bounds on individual
relaxations due to incomplete Benders iterations on a wider class of problems, the
numerical experiments presented show that the Branch-and-Benders approach promises
significant time savings for problems involving many scenarios.

In combination with more traditional parallelization avenues mentioned in the text, the
Branch-and-Benders approach clears the way for future high-performance implementations on
cluster computers.

\section{Acknowledgments}
\label{sec:acknowledgements}
This work was supported by the Swiss Commission for Technology and Innovation
(CTI), (Grant 16946.1 PFIW-IW). We thank the team at Adaptricity (Stephan Koch,
Andreas Ulbig, Francesco Ferrucci) for providing the system data and valuable
discussions on power systems.